\documentclass[11pt,reqno]{amsart}
\usepackage{diagrams}
\usepackage{url, cite}
\usepackage{amssymb}
\usepackage{enumerate}
\usepackage[pagewise]{lineno}
\numberwithin{equation}{section}

\theoremstyle{definition}
\newtheorem{thm}{Theorem}[section]
\newtheorem{cor}[thm]{Corollary}
\newtheorem{lem}[thm]{Lemma}

\newtheorem{prop}[thm]{Proposition}

\newtheorem{rem}[thm]{Remark}
\newtheorem{note}[thm]{Notation}

\DeclareMathOperator{\red}{\mathrm{red}}

\DeclareMathOperator{\Hc}{\mathcal{H}om}

\DeclareMathOperator{\p3}{\mathbb{P}^3}

\DeclareMathOperator{\pr}{\mathrm{pr}}

\DeclareMathOperator{\Spec}{\mathrm{Spec}}

\DeclareMathOperator{\I}{\mathcal{I}}
\DeclareMathOperator{\mo}{\mathcal{O}}

\newcommand{\mr}[1]{\mathrm{#1}}
\newcommand{\mb}[1]{\mathbb{#1}}
\newcommand{\mc}[1]{\mathcal{#1}}
\newcommand{\ov}[1]{\overline{#1}}

\begin{document}

\title{Determinant morphism for singular varieties}

\author[A. Dan]{Ananyo Dan}

\address{School of Mathematics and Statistics, University of Sheffield, Hicks building, Hounsfield Road, S3 7RH, UK}

\email{a.dan@sheffield.ac.uk}

\author[I. Kaur]{Inder Kaur}

\address{Institut für Mathematik,
Goethe-Universität Frankfurt,
Robert-Mayer-Str. 6-8,
60325 Frankfurt am Main, Germany}

\email{kaur@math.uni-frankfurt.de}

\subjclass[2010]{Primary $14$D$20$, $14$J$60$, $14$H$10$, Secondary $14$L$24$, $14$D$22$, $14$C$05$}


\date{\today}

\begin{abstract}
Let $X$ be a projective variety (possibly singular) over an algebraically closed field of any characteristic
and $\mc{F}$ be a coherent sheaf.
In this article, we define the determinant of $\mc{F}$ 
such that it agrees with the classical definition of 
determinant in the case when $X$ is non-singular. We study how the 
Hilbert polynomial of the determinant varies in families of 
singular varieties. Consider a singular family such 
that every fiber is a normal, projective variety.
Unlike in the case when the family is smooth,
the Hilbert polynomial of the determinant does not remain 
constant in singular families. 
However, we show that it exhibits
an upper semi-continuous behaviour. 
Using this we give a determinant morphism defined over flat families of 
coherent sheaves. This morphism coincides with  
the classical determinant morphism in the smooth case.
Finally, we give applications of our results to moduli spaces of semi-stable 
sheaves on $X$ and to Hilbert schemes of curves. 
\end{abstract}

\maketitle

\section{Introduction}
Let $X$ be a smooth, projective 
variety over an algebraically closed field of any characteristic. 
Classically, the determinant 
of a coherent sheaf $\mc{F}$ on $X$ is defined as 
$\otimes \det(F_i)^{(-1)^i}$, where $F_i$ is the $i$-th term in 
 a locally free resolution of $\mc{F}$ (see \cite[$1.1.17$]{huy}). Recall, when $X$ is a  
 smooth variety, the locally-free resolution of $\mc{F}$ is
 finite (see \cite[Proposition $2.1.10$]{huy}).
 Using this, one defines the determinant morphism 
 \[\det: M_X(P) \to \mr{Pic}(X),\]
 where $M_X(P)$ denotes the moduli space of (semi)stable sheaves
 with fixed Hilbert polynomial $P$
on $X$ and $\mr{Pic}(X)$ denotes the Picard scheme of $X$.
One application of the determinant morphism is that it 
can be used to define a moduli space of (semi)stable sheaves on 
$X$ with fixed
determinant arising as a fiber of the determinant morphism
$\det$ (see \cite[Theorem $4.5.4$]{huy}).
Note that, in characteristic $p$ this is possible only over the 
stable locus (see \cite[p. $381$]{langc}).
The moduli space of 
(semi)stable sheaves 
with fixed determinant often enjoys geometric properties not shared
by the space \emph{without} it. 
For example, in the case of $X$ a smooth, projective curve,
the moduli space of semi-stable vector bundles on $X$
with rank and degree coprime and fixed determinant is a Fano variety (see 
\cite{ind} for a proof). It is also rational (see \cite{ind2, king}). For $X$ a smooth, projective surface over a field of 
characteristic $0$ and $H$ an ample divisor on $X$,
the moduli space of rank $2$, semi-stable sheaves with fixed determinant 
$c_1$ and second Chern class $c_2 \gg 0$, denoted $M_{X,H}(2,c_1,c_2)$,
is irreducible (see 
\cite{barthmod} and \cite{giesli}). We also know that 
$M_{X,H}(r,c_1,c_2)$ is generically smooth for
$2rc_2-(r-1)c_1^2$ large enough (see \cite{donpoly, giesli2, ogmod}). Moreover, under certain assumptions, $M_{X,H}(2,c_1,c_2)$ is 
rational for the case $X=\mb{P}^2$ (see \cite{barthmod, hulek1, elstr}).
In \cite{langc}, the results on generic smoothness and irreducibility have 
been generalized to the case when the underlying field is 
of arbitrary characteristic.
Another frequent application
of determinant morphism is to study the birational 
geometry of $M_X(P)$ using the birational geometry 
of $\mr{Pic}(X)$ (see \cite{bayemmp, casablow}).  The geometric and topological properties of moduli spaces of semi-stable sheaves with 
fixed determinant on a nodal curve have been studied in \cite{mumf, genpre, hodc}.  

We know there exist moduli spaces of semi-stable sheaves
with fixed Hilbert polynomial in the case 
when $X$ is a singular, projective variety 
(see \cite{simp1, simp2, langa, langmix}).
For certain singular varieties of dimension $1$,
we also know the existence of semi-stable sheaves
with a prescribed determinant (see \cite{K4}).
It is natural to ask, whether one can define 
a determinant morphism for singular, projective varieties in a way that agrees with the definition of the
det morphism given above for smooth varieties.
Although a modification of the definition of determinant
for sheaves over an irreducible nodal curve $X$
has been given in \cite{bhos2} (see also 
\cite{nagsesh1}), this cannot be used to
define a determinant morphism on the entire moduli space  
(see \cite[Remark $4.8$]{bhos2}). This definition also does not extend to sheaves on varieties of higher dimension.
Another possible approach would be to study moduli spaces of principal $SL_r$ bundles. When $X$ is a smooth, projective curve, the moduli space of rank $r$ 
semi-stable sheaves with fixed determinant of degree coprime to $r$ is isomorphic to 
the moduli space of principal $\mr{SL}_r$ bundles.
However, this analogy does not extend to the case of nodal 
curves (see \cite[Theorem $3.2$ and Proposition $3.3$]{sct4}) 
(or even to the case of smooth, higher dimensional varieties).
Therefore, the question of constructing a moduli space of 
(semi)stable sheaves with fixed determinant 
over an arbitrary singular projective variety remains open.
In this article, we define a determinant morphism 
 when $X$ is a normal variety
of any dimension over a field of arbitrary characteristic
and give applications to geometry.

Let $k$ be an algebraically closed field of any characteristic and $X$ be 
a projective $k$-variety of arbitrary dimension (by variety we mean integral
of finite type over $k$).
Given a 
coherent sheaf $\mc{F}$ over $X$ of rank $n$, we define the \emph{determinant of }$\mc{F}$, denoted $\det(\mc{F})$, as the 
reflexive hull of $\wedge^n \mc{F}$ i.e., $\det(\mc{F}):= 
(\wedge^n \mc{F})^{\vee \vee}$. 
If $X$ is non-singular, then this definition of 
determinant agrees with the classical definition described earlier. Given two 
 Hilbert polynomials $L_1, L_2$, we say $L_1 \ge L_2$ if 
 $L_1(n) \ge L_2(n)$ for $n \gg 0$. Recall, if 
$\pi:\mc{X} \to S$ is a smooth family of projective varieties
with $S$ connected,
and $\mc{F}_S$ is an $S$-flat, coherent sheaf on $\mc{X}$,
then the Hilbert polynomial of the determinant of 
$\mc{F}_s:=\mc{F}_S|_{\mc{X}_s}$ remains fixed
as $s$ varies over the points in $S$. This 
follows from the existence of a finite, locally-free 
resolution of $S$-flat, coherent sheaves in smooth, 
projective families (see \cite[Proposition $2.1.10$]{huy}).
In contrast, when the family $\pi$ is not smooth, we 
prove the following:

\begin{thm}[Upper semi-continuity]\label{th:comp01}
 Let $\pi:\mc{X} \to S$ be a flat, projective morphism 
 such that for every $s \in S$, 
 the fiber $\mc{X}_s:=\pi^{-1}(s)$
 is integral, normal. Let $\mc{F}_S$ be a coherent sheaf on $\mc{X}$, flat over $S$ such that for every $s \in S$, the restriction 
 $\mc{F}_s$ of $\mc{F}_S$ to the fiber $\mc{X}_s$ is torsion-free. 
  Given a point $o \in S$, if the Hilbert 
 polynomial of $\det(\mc{F}_o)$ equals $L$, then 
 there exists an open neighbourhood $U \subset S$ containing
 $o$ such that for every $u \in U$, the Hilbert polynomial of 
 $\det(\mc{F}_u)$ is less than or equal to $L$.
\end{thm}

See Theorem \ref{c2} for the proof of a more general statement.
The main reason why the constancy of the Hilbert polynomial
fails in the setup of Theorem \ref{th:comp01} is that the 
dual (or double-dual) does not in general commute with 
pull-back. 
We prove that the upper semi-continuity of the 
Hilbert polynomials gives rise to a natural stratification 
of $S$ such that over each stratum the Hilbert polynomial 
of the determinant 
remains fixed. Moreover the stratification 
satisfies the following universal property: 
for any morphism from a 
noetherian scheme, say $S'$, to $S$, if the 
pull-back of $\mc{F}_S$ to $S'$
has constant Hilbert polynomial 
of the determinant over each fiber, then the morphism to $S$
must factor through this stratification. 
The theorem below is motivated by the classical flattening 
stratification, which has many applications
(see proofs of \cite[Theorem $2.2.4$ \mbox{ and }
Lemma $3.1.1$]{huy}), but has never so far been used to study the 
determinant morphism.
We 
show:

\begin{thm}[Stratification]\label{th:strintro}
Notations as in Theorem \ref{th:comp01}. 
 Given $s \in S$, denoted by $L(\mc{F}_s)$ the Hilbert 
 polynomial of the determinant of $\mc{F}_s$. Then, the set
 \[\mc{P}:=\{L(F_s)\, |\, s \in S\}\]
 is finite. Moreover, to each Hilbert polynomial $L \in \mc{P}$,
 one can associate a locally closed subscheme $S_L \subset S$
 satisfying the following conditions:
 \begin{enumerate}
  \item The natural morphism 
  $j:\coprod\limits_{L \in \mc{P}} S_L \to S$
  is bijective.
  \item For every $L \in \mc{P}$,  
  $\det(\mc{F}_{S_L})$ is flat over $S_L$ and for 
  any $s \in S_L$, $\det(\mc{F}_s)$ has Hilbert polynomial $L$
  and  is isomorphic to  $\det(\mc{F}_{S_L}) \otimes k(s)$.
  \item For any $L \in \mc{P}$, the determinant $\det(\mc{F}_{S_L})$
is $(S_2)$-relative to $\pi$ 
(see \S \ref{sec:ser} for the definition of $(S_2)$-relative to a morphism),
  \item If $g:S' \to S$ is a morphism from a reduced scheme 
  $S'$ to $S$, then 
  $g_{\mc{X}}^*\det(\mc{F}_S)$ is flat over $S'$,  where 
  $g_{\mc{X}}: \mc{X} \times_{S} S' \to \mc{X}$ is 
  induced  by $g$.
  Moreover, $g_{\mc{X}}^*\det(\mc{F}_S)$ is  $(S_2)$-relative to $\pi': \mc{X} \times_S S' \to S'$
    if and only if $g$ factors through $j$,
 \end{enumerate}
\end{thm}

See Theorem \ref{th:str} for a proof. 
Recall, that the determinant of a coherent sheaf is a rank $1$, torsion-free
sheaf. Therefore, it is a natural question to ask if there exists a determinant
morphism from $S$ to an appropriate
moduli space of rank $1$, torsion-free sheaves? Apriori, this 
has a negative answer because the definition of determinant
given above is not functorial. More precisely, given a morphism 
$\phi:T \to S$, the pull-back via $\phi$ of the determinant of 
$\mc{F}_S$ need \emph{not} be isomorphic to 
the determinant of the pull-back
of $\mc{F}_S$ via $\phi$. However,
as a corollary to Theorem \ref{th:strintro}, we can show the following. Let $\ov{\mr{Pic}}_{\pi,L}$ be the relative compactified Picard scheme parametrizing
rank $1$, torsion-free
sheaves on fibers of $\mc{X}$ with Hilbert polynomial $L$ (see \S \ref{sec:pic}). In Corollary \ref{det11}, we show that there exists  an \emph{unique}
determinant morphism from the stratification of $S$
(mentioned in Theorem \ref{th:strintro}) to the disjoint
union of relative compactified Picard schemes $\ov{\mr{Pic}}_{\pi,L}$
(as $L$ varies over Hilbert polynomials in $\mc{P}$),
satisfying the 
obvious universal property.

{\bf{Application to moduli spaces}:}
We use the determinant morphism given in 
Theorem \ref{th:strintro} to define a 
moduli space of semi-stable sheaves with fixed determinant
over normal, projective varieties.
Let $X$ be a normal, projective variety, $P$ a 
Hilbert polynomial 
of a rank $n$ torsion-free sheaf on $X$ with degree coprime
to $n$. By \cite{langmix, simp1, simp2} there exists a projective scheme of finite type parametrizing such sheaves. Denote this scheme by $\mc{M}_X(P)$. Let $\mc{L}$ be a 
rank $1$, torsion-free sheaf on $X$. We define 
a moduli functor with fixed determinant $\mc{L}$, denoted 
$\mc{M}_X(P,\mc{L})$, which to a scheme $S$ over $k$, associates
the set of $S$-flat, semi-stable sheaves $\mc{F}_S$ over 
$X \times_k S$ with fiber-wise Hilbert polynomial $P$ and 
determinant $\mc{L}$ and $\det(\mc{F}_S)$ is isomorphic to 
$\mc{L} \otimes_k \mo_S$. The last condition 
ensures functoriality (see Lemma  \ref{lem:func}).
By Theorem \ref{th:strintro} above, there is a 
natural stratification
\[\coprod_{L \in \mc{P}} M_X(P,L)\]
of $M_X(P)_{\red}$ by the set
$\mc{P}$. Here $M_X(P,L)$ denotes the subscheme 
parametrizing semi-stable sheaves $\mc{F}$ on $X$ with Hilbert 
polynomial $P$ and the Hilbert polynomial of $\det(\mc{F})$
equals $L$. 
The universal property of the compactified 
Picard functor gives rise to the determinant morphism 
\[\mr{det}_{M_X(P)_{\red}}: 
\coprod_{L \in \mc{P}} {M}_X(P,L) \to \coprod_{L
\in \mc{P}} \ov{\mr{Pic}}_{X,L}.\]
The natural candidate for corepresenting 
$\mc{M}_X(P,\mc{L})$ is the fiber, denoted by $M_X(P,\mc{L})$, to the determinant morphism $\det_{M_X(P)_{\red}}$
 over the point 
corresponding to $\mc{L}$.
 We  prove:

\begin{thm}[See Theorem \ref{c13}]\label{cn25}
 Let $L$ be the Hilbert polynomial of $\mc{L}$. We have the following:
 \begin{enumerate}
  \item Suppose $\mc{M}_X(P,\mc{L})$ is corepresentable by a locally closed subscheme of $M_X(P)$, which we denote by $N_{\mc{L}}$. Then 
 $N_{\mc{L},\red}=M_X(P,\mc{L})_{\red}$.
 \item If $M_X(P)$ is reduced and $L$ is the smallest 
 Hilbert polynomial in the set $\mc{P}$, then 
 $\mc{M}_X(P,\mc{L})$ is 
 universally corepresentable by the quasi-projective 
 $k$-scheme  $M_X(P,\mc{L})$.
 \end{enumerate}
\end{thm}

{\bf{Application to geometry}:}
In the seminal work by Schlessinger and Piene 
\cite{piesch} (see also \cite{elli3}), the authors show that the Hilbert 
scheme (of closed subschemes in $\p3$)
corresponding to the Hilbert polynomial
$P(n)=3n+1$ consists of two irreducible components,
one parameterizing twisted cubic curves and
the other parameterizing plane cubic curves 
along with an embedded point.
Several aspects of this Hilbert scheme
have been studied (see \cite{elli4} and more recently 
\cite{bxia}).
Using Theorem \ref{th:strintro} 
(replacing $\mc{F}_S$ by the 
ideal sheaf of the family of curves), we show that this
Hilbert scheme has a non-trivial stratification.
We can check that the twisted cubic curves are Cohen-Macaulay
but the plane cubic curves with an embedded point are not
(see \S \ref{ex:sch} for further details). 
Motivated by this we ask the following general question:
given a flat family of projective curves,  
how ``far'' are the fibers from being Cohen-Macaulay?
By a \emph{curve} we mean a closed subscheme of dimension $1$. It can be non-reduced or have embedded points.

To study this we introduce the notion of 
Cohen-Macaulayfication.
 Given a normal, projective,
integral surface $X$ and a curve $C \subset X$, 
we define the \emph{Cohen-Macaulayfication of} $C$ to be 
\[\mo_C^{\mr{CM}}:=\mc{E}xt^1_X(\mc{E}xt^1_X(\mo_C,\mo_X),
 \mo_X).\]
 Note that, $\mo_C^{\mr{CM}}$ is a Cohen-Macaulay 
 $\mo_X$-module.
If $C$ is Cohen-Macaulay to begin with, then 
$\mo_C \cong \mo_C^{\mr{CM}}$. Moreover, if $C$ is 
generically reduced, then $\mo_C^{\mr{CM}}$ does not 
depend on the choice of $X$.
Applying Theorem \ref{th:strintro} to Hilbert schemes of curves contained in normal, integral, projective
surfaces, we obtain a necessary
and sufficient condition for a fiber to be 
Cohen-Macaulay. 
In particular,
if all the fibers were Cohen-Macaulay, there would 
be no non-trivial stratification. More precisely, 
we show:

\begin{thm}[see Theorem \ref{th:geom}]\label{th:cmintro}
Fix Hilbert polynomials $P$ and $Q$ of degrees $1$ and 
$2$, respectively. 
 Let $(\mc{C}_S \subset \mc{X}_S \subset \mb{P}^3_S)$
 be a flat family of closed subschemes of $\p3$ such that
 for any $s \in S$, the triple restricts to the 
 inclusions $(\mc{C}_s \subset \mc{X}_s \subset 
 \mb{P}^3_s)$, where $\mc{C}_s$ (resp. $\mc{X}_s$)
 is closed subscheme of $\mb{P}^3_s$ with Hilbert 
 polynomial $P$ (resp. $Q$) and $\mc{X}_s$ is normal
 and integral. Denote by $\I_{\mc{C}_S}$ the ideal sheaf
 of $\mc{C}_S$ as a closed subscheme of $\mc{X}_S$. Let $\mc{P}$
 be as in Theorem \ref{th:strintro} and substitute
 $\mc{F}_S$ by $\I_{\mc{C}_S}$. Let $\{S_L\}_{L \in \mc{P}}$
 be the associated stratification of $S$. Then,
 \begin{enumerate}
  \item for every $L \in \mc{P}$, we have $L \ge Q-P$
  with equality if and only if there exists $s \in S_L$ 
  for which $\mc{C}_s$ is Cohen-Macaulay,
  \item $s \in S_L$ if and only if the 
  Cohen-Macaulayfication $\mo_{\mc{C}_s}^{\mr{CM}}$ of 
  $\mc{C}_s$ has Hilbert polynomial $Q-L$.
 \end{enumerate}
\end{thm}

\vspace{0.2cm}
\emph{Notations}:
For $S$-schemes $\mc{X}_S$, $T$ denote by $\mc{X}_T:=\mc{X}_S \times_S T$. For any sheaf $\mc{F}_S$ on $\mc{X}_S$, denote by $\mc{F}_T:=\mc{F}_S \otimes_S \mo_T$.
For any morphism of sheaves $f_S:\mc{F}_S \to \mc{G}_S$, denote by $f_T:\mc{F}_T \to \mc{G}_T$ the natural morphism $f_S \otimes \mr{id}$.

\vspace{0.2 cm}
\emph{Acknowledgements}:
We thank Prof. A. Langer for his feedback and 
Prof. C. Simpson and Prof. D. S. Nagaraj
for helpful discussions. 
The first author is funded by EPSRC grant number EP/T019379/1 and the second author
by the DFG, TRR $326$
Geometry and Arithmetic of Uniformized Structures, project number $444845124$.
A part of the work was done when the first author was supported by ERCEA Consolidator Grant $615655$-NMST and also
by the Basque Government through the BERC $2014-2017$ program and by Spanish
Ministry of Economy and Competitiveness MINECO: BCAM Severo Ochoa
excellence accreditation SEV-$2013-0323$ and the second author by a PNPD fellowship from CAPES, Brazil.

\section{Determinant of a coherent sheaf on a singular variety}
In this section we introduce the notion of determinant of 
coherent sheaves over singular varieties. 

\subsection{Relative $(S_k)$-sheaves}\label{sec:ser}
Recall, Serre's criterion $(S_k)$: Given a noetherian scheme
$Y$ and a coherent sheaf $\mc{G}$ on $Y$, we say that 
$\mc{G}$ \emph{has property} $(S_k)$ if for every $y \in Y$, 
 \[\mr{depth}_{\mo_{Y,y}}(\mc{G}_y) \ge 
     \min \{k,\dim(\mr{supp}(\mc{G}_y))\}.     \]
If there exists a flat, 
projective morphism $f:Y \to S$, we say that an $S$-flat,
coherent sheaf $\mc{G}$ on $Y$, is 
$(S_k)$-\emph{relative to} $f$
if for every $y \in Y$, $s=f(y)$, we have 
\[\mr{depth}_{\mo_{Y,y}}((\mc{G}|_{Y_s})_y) \ge \min\{k,
 \dim \mo_{Y,y}-\dim \mo_{S,s}\}.\]
See \cite[\S $3$]{hass} for interesting properties of such 
sheaves.

Note that in the case when every fiber of $f$ is normal, integral
and $\mc{G}$ is $(S_1)$-relative to $f$, there exists a 
closed subscheme $Z \subset Y$ such that for every 
$s \in S$, $\mr{codim}(Z \cap Y_s,Y_s)\ge 2$ and 
$\mc{G}|_{Y\backslash Z}$ is locally-free. Indeed, for any 
$s \in S$, the restriction $\mc{G}_s$ of $\mc{G}$ to $Y_s$ is 
torsion-free. As $\mc{G}_s$ is torsion-free, there exists 
an open subset $U_s \subset Y_s$ which is a complement 
of a subvariety of codimension at least $2$ and 
$\mc{G}_s|_{U_s}$
is locally-free. Using \cite[Theorem $14.23$]{gortz}, we 
conclude that $\mc{G}$ is locally-free at every point 
of $U_s$. As local freeness is an open condition, there exists an 
open subscheme $\widetilde{U}_s \subset Y$ containing $U_s$
such that $\mc{G}|_{\widetilde{U}_s}$ is locally-free.
Since $s \in S$ is arbitrary, taking the union $U$ of the 
open subschemes 
$\widetilde{U}_s$ as $s$ varies over points on $S$,
we observe that $\mc{G}|_U$ is locally-free. Moreover, for every 
$s \in S$, the complement $Y_s\backslash (U \cap Y_s)$
is of codimension at least $2$ in $Y_s$. We denote by 
$Z:=Y \backslash U$ the complement of $U$ in $Y$ and we say that 
$Z$ \emph{is of fiberwise codimension at least} $2$.

\subsection{Setup}\label{se:det}
Let $A$ be a noetherian $k$-algebra ($k$ is an algebraically closed field
of any characteristic) and $S:=\Spec(A)$. Suppose that $S$ is integral.
 Let  
\[\delta:\mc{X} \to S\]
be a flat, projective morphism such that for 
every $s \in S$, the fiber $\mc{X}_s:=\delta^{-1}(s)$ is normal and integral. 
Fix a coherent sheaf $\mc{F}_S$ on $\mc{X}$, flat over $S$ such that for every $s \in S$, the restriction $\mc{F}_s:=\mc{F}_S|_{\mc{X}_s}$ is torsion-free and of rank $n$.
In other words $\mc{F}_S$ is $(S_1)$-relative to $\delta$. 

\subsection{Determinant of a coherent sheaf} 
Let $X$ be a projective scheme  and $\mc{F}$ be a 
coherent sheaf on $X$. For any $n \ge 0$, denote by 
\[\mr{det}_n(\mc{F}):= 
(\wedge^n \mc{F})^{\vee \vee}.\]
When $X$ is integral, 
 of finite type over a field, then the rank of a coherent sheaf on $X$
 is well-defined. In this case, if $\mc{F}$ is of rank $n$, we 
 call $\det_n(\mc{F})$ the \emph{determinant of} $\mc{F}$ and 
 denote it simply by $\det(\mc{F})$. Since the dual of a coherent sheaf is torsion-free (see \cite[Corollary $1.2$]{stabhar}),  
   $\det(\mc{F})$ is torsion-free of rank $1$. 
Since rank one, torsion-free sheaves on projective varieties are semi-stable, it is semi-stable. 
Note that, if $\mc{F}$ is a locally-free sheaf or is of finite 
projective dimension, then $\det(\mc{F})$ is the same as the classical 
definition of determinant as mentioned in the Introduction.

\subsection{Functoriality of determinant}
Let $S$ be a noetherian scheme and $X$ be a normal, projective variety. Denote by $\mc{F}_S$ a coherent sheaf on $X \times_k S$ of rank $n$, flat over $S$ and 
$(S_1)$-relative to the natural projection $f:X \times S
\to S$.
We claim that $\det_n(\mc{F}_S)$ is reflexive.
Indeed, consider the short exact sequence:
\begin{equation}\label{eqn14}
 0 \to (\wedge^n \mc{F}_S)/\mr{tors} \to (\wedge^n 
 \mc{F}_S)^{\vee \vee} \to T \to 0.
\end{equation}
Using \S \ref{sec:ser}, note that $T$ is supported 
on a closed subscheme $Z \subset X \times_k S$ of 
fiberwise codimension at least $2$. By 
\cite[Proposition $1.2.3$ and Lemma $1.2.4$]{brun1}, this 
implies \[\Hc_{X \times S} (T, \mo_{X \times S})=0=
    \mc{E}xt^1_{X \times S}(T,\mo_{X \times S}).\]
Then, dualizing \eqref{eqn14}, we get 
$(\wedge^n \mc{F}_S/\mr{tors})^\vee$ is isomorphic to 
$\det_n(\mc{F}_S)^\vee$. Moreover,  observe that 
$(\wedge^n \mc{F}_S/\mr{tors})^\vee$ is isomorphic to 
$(\wedge^n \mc{F}_S)^\vee$. This implies,
\[(\wedge^n \mc{F}_S/\mr{tors})^{\vee \vee} \cong 
 (\wedge^n \mc{F}_S)^\vee \cong \mr{det}_n(\mc{F}_S)^{\vee
 \vee}.\]
This proves the claim. We now use this to show:

  \begin{lem}\label{lem:func}
  Let $S:=\Spec(A)$ be a noetherian, affine scheme
  and $f:X \times_k S \to S$ be the natural projection,
   where $X$ is a normal, projective variety. Assume that $\mc{F}_S$ is as above and  
   for every
  $s \in S$, the 
  natural morphism \[\mr{det}_n(\mc{F}_S) \otimes_{\mo_S} k(s)
    \to \det(\mc{F}_s)\]
    is an isomorphism. Then, for any ring homomorphism
    (of finite-type) $A \to B$,
    we have  \[\phi_B: \mr{det}_n(\mc{F}_S) \otimes_{\mo_S} B
    \to \mr{det}_n(\mc{F}_S \otimes_{\mo_S} B)\]
    is an isomorphism.
  \end{lem}

  \begin{proof}
    Using \S \ref{sec:ser}, there exists a closed 
    subscheme $Z \subset X \times_k S$ of fiberwise
    codimension at least $2$ and $\mc{F}_S|_{X \times
    S \backslash Z}$ is locally-free. 
    Denote by $U:=X \times_k S \backslash Z$ and 
    $U_B$ the pull-back of $U$ under the natural morphism from 
    $X \times \Spec(B)$ to $X \times \Spec(A)$.
    Since duality of locally-free sheaves commute with 
    pull-back, the morphism $\phi_B$ is an isomorphism over
    $U_B$. Since $\det(\mc{F}_s)$ is reflexive 
    for any $s \in S$, we have by assumption that 
    $\det_n(\mc{F}_S)$ is $(S_2)$-relative to $f$. 
    In particular, $\det_n(\mc{F}_S) \otimes_A B$ is 
    $(S_2)$-relative to $f_B: X \times_k \Spec(B) \to 
    \Spec(B)$. We observed above that 
    $\det_n(\mc{F}_S \otimes_A B)$ is reflexive.
    Then, \cite[Proposition $3.6$]{hass} implies that $\phi_B$
    is an isomorphism. This proves the lemma.
           \end{proof}

\section{Determinants of sheaves in families}
In this section we prove that the Hilbert polynomial of 
the determinant of coherent sheaves in flat families of singular varieties satisfies the 
upper semi-continuity property (Theorem \ref{c2}).

\begin{lem}\label{c1}
Suppose that $A$ in the setup is a discrete valuation ring. Denote by $k$ the residue field of $A$.
 Let $\mc{G}$ be a torsion-free sheaf on $\mc{X}$. Then  the natural morphism 
 \[\mc{G}^\vee \otimes_A k \to (\mc{G}|_{\mc{X}_k})^\vee\]
 is injective, where $\mc{X}_k$ is the special fiber of $\delta$.
\end{lem}

\begin{proof}
 Since $\mc{G}$ is a coherent sheaf, there exists $n \gg 0$ such that $\mc{G}(n)$ is globally generated (see \cite[Theorem II.$5.17$]{R1}) i.e.,
 there exists an indexed set $I$ of finite cardinality such that we have a surjective morphism of the form 
 \[\mc{H}_0:=\bigoplus\limits_{i \in I} \mo_{\mc{X}} \twoheadrightarrow \mc{G}(n). \]
 Consider the resulting short exact sequence:
 \[0 \to \mc{G}' \to \mc{H}_0(-n) \to \mc{G} \to 0\]
 for a coherent sheaf $\mc{G}'$ on $\mc{X}$.
 Taking dual, we get the short exact sequence:
 \[0 \to \mc{G}^\vee \to \mc{H}_0^\vee(n) \to \mc{G}'^\vee \xrightarrow{\phi} \mc{E}xt^1(\mc{G},\mo_{\mc{X}}) \to 0, \, \mbox{ where } (-)^\vee=\Hc_{\mc{X}}(-,\mo_{\mc{X}}).\]
 Since $\mc{G}'^\vee$ is torsion-free and $\ker \phi \subset \mc{G}'^\vee$, we have that $\ker \phi$ is torsion-free, hence $R$-flat.
 This implies $\mc{G}^\vee \otimes_R k \to \mc{H}_0^\vee(n) \otimes_R k$ is injective. We then have the following commutative diagram:
 \[\begin{diagram}
 \mc{G}^\vee \otimes_R k&\rInto&\mc{H}_0^\vee(n)\otimes_R k\\
  \dTo_{f_1}&\circlearrowleft&\dTo_{f_2}\\
    (\mc{G}|_{\mc{X}_k})^\vee&\rInto&(\mc{H}_0(-n)|_{\mc{X}_k})^\vee
       \end{diagram}\]
 where the injectivity of the lowest row follows from the left exactness of $\Hc(-,\mo_{\mc{X}_k})$ applied to the above short exact sequence.
Since $\mc{H}_0^\vee(n)$ is locally free, $f_2$ is an isomorphism.
By the commutativity of diagram, we then conclude that $f_1$ is injective.
This proves the lemma.
\end{proof}

Given two Hilbert polynomials $L_1$ and $L_2$, we say that 
$L_1 \ge L_2$ if for $n \gg 0$, we have $L_1(n) \ge L_2(n)$.
We prove:

\begin{thm}[Upper semi-continuity]\label{c2}
Denote by $L$ the Hilbert polynomial of the determinant of $\mc{F}_{\eta}$,
where $\eta$ is the generic point of $S$ and $\mc{F}_{\eta}$ is the 
restriction of $\mc{F}_S$ to the fiber $\mc{X}_\eta:=\pi^{-1}(\eta)$. Then, the following holds true:
\begin{enumerate}
\item for every $s \in S$, the Hilbert polynomial of the 
determinant of $\mc{F}_s$ is greater than or equal to $L$,
\item if $\det_n(\mc{F}_S) \otimes_{\mo_S} k(s)$ is reflexive, 
then $\det(\mc{F}_s)$ is isomorphic to $\det_n(\mc{F}_S) \otimes_S
k(s)$ and  is of Hilbert polynomial $L$.
\end{enumerate}
\end{thm}

\begin{proof}
Consider the natural morphism \[\phi_s: (\wedge^n \mc{F}_S)^\vee 
        \otimes_{\mo_S} k(s) \to (\wedge^n \mc{F}_s)^\vee. \]
Note that, there exists a closed subscheme $Z \subset \mc{X}$
of fiberwise codimension at least $2$ as 
mentioned in \S \ref{sec:ser} such that 
$\mc{F}_S|_{\mc{X}\backslash Z}$ is locally-free.
Denote by $U:=\mc{X} \backslash Z$ and $U_s:=U \cap \mc{X}_s$.
As duality for locally-free sheaves 
commute with pull-back, the morphism $\phi_s$ is an isomorphism
over $U_s$. Dualizing $\phi_s$, we get a morphism from  
 $\det(\mc{F}_s)$ to
$((\wedge^n \mc{F}_S)^\vee \otimes_{\mo_S} k(s))^\vee$
which is an isomorphism over $U_s$.
Since both sheaves are reflexive, \cite[Proposition $1.6$]{stabhar}
implies that 
\begin{equation}\label{eqn10}
 ((\wedge^n \mc{F}_S)^\vee \otimes_{\mo_S} k(s))^\vee \cong 
\det(\mc{F}_s).
\end{equation}
Consider now a discrete valuation ring $R_s$ 
(depending on the point $s$) 
along with a morphism $\mr{sp}_s: \Spec(R_s)
\to S$ such that the generic point (resp. special point) of 
$\Spec(R_s)$ maps to $\eta$ (resp. $s$). Denote by $\mc{F}_{R_s}$ 
the pull-back of $\mc{F}_S$ to $\mc{X}_{R_s}$ via 
the specialization morphism $\mr{sp}_s$.
Note that, $\mc{F}_{R_s}$ is flat over $\Spec(R_s)$.
By Lemma \ref{c1}, we have an injective morphism 
\begin{equation}\label{eqn11}
 (\wedge^n \mc{F}_{R_s})^{\vee \vee} \otimes_{R_s} k(s) 
 \hookrightarrow ((\wedge^n \mc{F}_{R_s})^\vee \otimes_{R_s} 
 k(s))^\vee
\end{equation}
Using \eqref{eqn10} (after substituting $S$ by $\Spec(R_s)$)
we have $((\wedge^n \mc{F}_{R_s})^\vee \otimes_{R_s} 
 k(s))^\vee$ is 
 isomorphic to $\det(\mc{F}_s)$. As $\det_n(\mc{F}_{R_s})$ is 
 $R_s$-flat, the Hilbert polynomial of $\det_n(\mc{F}_{R_s})
 \otimes_{R_s} k(s)$ is the same as that of 
 $\det_n(\mc{F}_{R_s}) \otimes_{R_s} K$, where $K$ (resp. $k(s)$)
 is the fraction field (resp. residue field) of $R_s$. 
 As duality commutes with flat base change, we have that
 $\det_n(\mc{F}_{R_s}) \otimes_{R_s} K$ is isomorphic to 
 $\det_n(\mc{F}_K)$. By assumption, $\mc{F}_K \cong \mc{F}_\eta$,
 the determinant of which is of Hilbert polynomial $L$.
 Hence, \eqref{eqn11} implies that the Hilbert polynomial of 
 $\det(\mc{F}_s)$ is greater than or equal to $L$. This proves 
 $(1)$.
 
Consider the natural morphism 
\begin{equation}\label{eqn12}
 \mr{det}_n(\mc{F}_S) 
        \otimes_{\mo_S} R_s \to ((\wedge^n \mc{F}_S)^\vee 
        \otimes_{\mo_S} R_s)^\vee
\end{equation}
 Let 
$U_{R_s} \subset \mc{X}_{R_s}$ be the pull-back of $U$ under the
specialization morphism $\mr{sp}_s$ from $\Spec(R_s)$ to $S$.
Since duality of locally-free sheaves commute with pull-back,
the natural morphism 
\[\phi_{R_s}: (\wedge^n \mc{F}_S)^\vee \otimes_{\mo_S}
 R_s \to (\wedge^n \mc{F}_{R_s})^\vee\]
 is an isomorphism over $U_{R_s}$.
 Dualizing $\phi_{R_s}$, we get a morphism from 
 $\det_n(\mc{F}_{R_s})$ to 
 $((\wedge^n \mc{F}_S)^\vee \otimes_{\mo_S} R_s)^\vee$, 
 which is an isomorphism over $U_{R_s}$.
 Since the complement of $U_{R_s}$ in $\mc{X}_{R_s}$ is of 
 codimension at least $2$ and both sheaves are reflexive, 
  \cite[Proposition $1.6$]{stabhar} gives an isomorphism 
\[\mr{det}_n(\mc{F}_{R_s}) \xrightarrow{\sim} 
 ((\wedge^n \mc{F}_S)^\vee \otimes_{\mo_S} R_s)^\vee.\]
Substituting this in \eqref{eqn12} and applying $- \otimes_{R_s}
k(s)$, we get the following composed morphism 
\[\mr{det}_n(\mc{F}_S) \otimes_{\mo_S} k(s) \cong 
\mr{det}_n(\mc{F}_S) \otimes_{\mo_S} R_s \otimes_{R_s} k(s)
 \xrightarrow{\phi_1} \mr{det}_n(\mc{F}_{R_s}) \otimes_{R_s}
 k(s) \xrightarrow{\phi_2} \det(\mc{F}_s),\]
where $\phi_2$ is injective as observed above using 
\eqref{eqn10} and \eqref{eqn11}.
Note that, $\det_n(\mc{F}_S) \otimes_S k(s)$, 
$\det_n(\mc{F}_{R_s}) \otimes_{R_s} k(s)$ and $\det(\mc{F}_s)$
agree over $U_s$.
Furthermore, $\det(\mc{F}_s)$ is reflexive. As a result, if 
$\det_n(\mc{F}_S) \otimes_{\mo_S} k(s)$ is reflexive, then 
by \cite[Proposition $1.6$]{stabhar} we have 
$\det_n(\mc{F}_S) \otimes_{\mo_S} k(s) \cong \det(\mc{F}_s)$.
 This means that the morphism $\phi_1$ is injective and 
 $\phi_2$ is surjective. As $\phi_2$ is already injective, 
 this means that $\phi_2$ is an isomorphism. Moreover,
 we observed above that the Hilbert polynomial of 
 $\mr{det}_n(\mc{F}_{R_s}) \otimes_{R_s} k(s)$ is $L$. This proves
 $(2)$, hence the theorem.
\end{proof}

 \section{The determinant morphism}
 In this section we use the upper semi-continuity result to give a natural stratification (satisfying an 
 universal property) of the parameterizing space $S$ of a flat 
 family of coherent sheaves such that the Hilbert 
 polynomial of the determinant remains constant 
 in each stratum (Theorem \ref{th:str}). 
 In \S \ref{sec:pic} we recall known facts about the relative 
 compactified Picard functor.
 We then give a determinant morphism 
 from the stratification of $S$
 to a union of relative compactified Picard schemes (Corollary \ref{det11}).
 
  \subsection{Setup}
Throughout this section, denote by $k$ an algebraically 
closed field of any characteristic. Let $S$ be 
a connected, reduced scheme
and \[\delta: \mc{X} \to S\] 
be a flat, projective morphism such that for 
every $s \in S$, the fiber 
$\mc{X}_s:=\delta^{-1}(s)$ is normal and integral. 
Fix a coherent sheaf $\mc{F}_S$ over $\mc{X}$, flat over $S$.
Assume that the restriction $\mc{F}_s$ of $\mc{F}_S$ to any fiber $\mc{X}_s:=
\delta^{-1}(s)$ is a torsion-free $\mo_{\mc{X}_s}$-module. 
Suppose that the rank of $\mc{F}_s$
is $n$, for any $s \in S$.

\subsection{Stratification by Serre's criterion}
We now show that the locus $V \subseteq S$, of  
points $s$ such that the Hilbert 
polynomial of $\det(\mc{F}_s)$ is minimal among all possible Hilbert 
polynomials of the determinant of the fibers of $\mc{F}_S$
coincides with the locus of points where 
$\det(\mc{F}_S)$ is $(S_2)$-relative to $\delta$. 
This ensures the openness
of $V$.
The complement of open sets are proper, 
closed sets.
We can recursively then apply the upper semi-continuity property
to this closed sub-scheme. Since $S$ is noetherian, we can check that 
this process terminates after finite steps.

\begin{lem}\label{lem:detconst}
 Let $V \subset S$ be a connected open subset such that $\det_n(\mc{F}_S)|_{\mc{X}_V}$ is $(S_2)$-relative to $\delta$ (restricted to 
$\mc{X}_V$). Then, the Hilbert polynomial of $\det(\mc{F}_v)$ remains
fixed as $v$ varies over the points on $V$.
\end{lem}

\begin{proof}
Denote by $\mc{F}_V:=\mc{F}|_{\mc{X}_V}$. Let $W$ be an 
 irreducible component of $V$. Consider
the natural morphism 
\[\phi_W: \mr{det}_n(\mc{F}_V)|_{\mc{X}_{W}}
    \to ((\wedge^n \mc{F}_V)^\vee \otimes \mo_{W})^{\vee}.\]
By \S \ref{sec:ser}, there exists a closed subscheme $Z \subset
\mc{X}$ of fiberwise codimension at least $2$, such that 
$\mc{F}_S|_{\mc{X}\backslash Z}$ is locally-free.
Denote by $U:=\mc{X}\backslash Z$. As duality for locally-free
sheaves commute with pull-back, the morphism $\phi_W$ is an 
isomorphism over $U_W:=U \cap \mc{X}_W$. 
Observe that $\mr{det}_n(\mc{F}_V)|_{\mc{X}_{W}}$ is 
$(S_2)$-relative
to $\delta$ (restricted to $\mc{X}_W$)
and $((\wedge^n \mc{F}_V)^\vee \otimes \mo_{W})^{\vee}$ is reflexive. By \cite[Proposition $3.6$]{hass} implies that 
$\phi_W$ is an isomorphism. Moreover, the natural morphism 
from $(\wedge^n \mc{F}_V)^\vee \otimes \mo_{W}$ to 
$(\wedge^n \mc{F}_W)^\vee$ is an isomorphism over $U_W$.
Then again by \cite[Proposition $3.6$]{hass}, dualizing this morphism gives us an isomorphism:
\[\mr{det}_n(\mc{F}_W) \xrightarrow{\sim}
 ((\wedge^n \mc{F}_V)^\vee \otimes \mo_W)^\vee.\]
Therefore, $\det_n(\mc{F}_W)$ is isomorphic to $\det_n(\mc{F}_V)|_{\mc{X}_W}$.
Hence $\det_n(\mc{F}_W)$ is $(S_2)$-relative to $\delta$ (restricted
to $\mc{X}_W$).
By Theorem \ref{c2}, for any $w \in W$,
the Hilbert polynomial of $\det(\mc{F}_w)$ is the same 
as the Hilbert polynomial of $\det(\mc{F}_\eta)$ for 
$\eta \in W$ the generic point. This proves the lemma.
\end{proof}

\begin{thm}[Stratification]\label{th:str}
 For $s \in S$, denote by $L(\mc{F}_s)$ the Hilbert 
 polynomial of the determinant of $\mc{F}_s$. Then, the set
 \[\mc{P}:=\{L(F_s)\, |\, s \in S\}\]
 is finite. Moreover, to each Hilbert polynomial $L \in \mc{P}$,
 one can associate a locally closed subscheme $S_L \subset S$
 satisfying the following conditions:
 \begin{enumerate}
  \item The natural morphism 
  $j:\coprod\limits_{L \in \mc{P}} S_L \to S$
  is bijective.
  \item For every $L \in \mc{P}$,  
  $\det_n(\mc{F}_{S_L})$ is flat over $S_L$, $(S_2)$-relative
  to $\delta$ (restricted to $\mc{X}_{S_L}$) and for 
  any $s \in S_L$, $\det_n(\mc{F}_{S_L}) \otimes k(s)$ is 
  isomorphic to $\det(\mc{F}_s)$ with Hilbert polynomial $L$.
  \item Given a morphism $g:S' \to S$ from a reduced scheme 
  $S'$ to $S$, we have that
  $g_{\mc{X}}^*\det_n(\mc{F}_S)$ is flat over $S'$,  where $g_{\mc{X}}: \mc{X} \times_{S} S' \to \mc{X}$ is 
  induced by $g$.
  Furthermore, $g_{\mc{X}}^*\det_n(\mc{F}_S)$ is $(S_2)$-relative to $\delta': \mc{X} \times_S S' \to S'$
  if and only if $g$ factors through $j$.
 \end{enumerate}
\end{thm}

\begin{proof}
 By the theorem of generic flatness (see \cite[Theorem $10.84$]{gortz}),
there exists an open dense subscheme $U \subset S$ such that the restriction 
$\det_n(\mc{F}_S)|_U$ of $\det_n(\mc{F}_S)$ to $\mc{X}_U:=
\delta^{-1}(U)$ is flat over 
$U$. Take $U$ to be the largest such open subset of $S$.
As duality commutes with flat pull-back, 
for any generic point $\eta \in U$, the natural morphism 
\[\mr{det}_n(\mc{F}_S)|_{\mc{X}_U} 
\otimes_{\mo_U} k(\eta) \to \mr{det}_n(\mc{F}_{\eta})\]
is an isomorphism. 
In particular, $\mr{det}_n(\mc{F}_S)|_{\mc{X}_U} 
\otimes_{\mo_U} k(\eta)$ is an $(S_2)$-sheaf. Using 
\cite[Theorem $12.2.1$]{ega43}, we conclude that there exists an 
open dense subscheme $V \subset U$ such that the restriction
$\det_n(\mc{F}_S)|_{\mc{X}_V}$ of $\det_n(\mc{F}_S)$ to 
$\mc{X}_V$ is $(S_2)$-relative to $\delta$ (restricted to 
$\mc{X}_V$). We take the largest 
such open subset $V$ of $U$. 
 Denote by
$T_1:= V$ and $\{L_{1,j}\}$ the set of Hilbert polynomials of 
$\det(\mc{F}_v)$ as $v$ varies over points in $V$. 
By Lemma \ref{lem:detconst}, the 
Hilbert polynomial of the determinant remains fixed over 
any connected component of $V$. Denote by $T_1^c$ the 
complement $S \backslash T_1$. Note that, $T_1^c$ is 
assumed to be reduced, by construction.
We continue recursively as follows: given a closed subscheme
$T_i^c$ of $S$, denote by $T_{i+1}$, the largest open subscheme
of $T_i^c$ such that 
$\det_n(\mc{F}_{T_i^c})|_{\mc{X}_{T_{i+1}}}$
is $T_{i+1}$-flat and $(S_2)$-relative to $\delta$ 
(restricted to $\mc{X}_{T_{i+1}}$). Denote by $T_{i+1}^c:=
T_i^c \backslash T_{i+1}$. As observed in the case $i=1$,
there are only finitely many Hilbert polynomials of the 
determinant of $\mc{F}_t$ as $t$ ranges over points in $T_i$
(this follows from the existence of finitely many connected 
components of $T_1$). Similarly, to each $i$, we 
associate the set of Hilbert polynomials $\{L_{i,j}\}$
of the determinant of $\mc{F}_t$ as $t$ varies over points 
in $T_i$. As $S$ is noetherian, there exists a finite $N$
for which $T^c_N$ coincides with $T_N$ i.e., $T_{N+1}^c=
\emptyset$ (use finiteness of descending chain of closed
subschemes). Hence, the set 
 \[\mc{P}:=\{L(F_s)\, |\, s \in S\}=\{L_{i,j}\}\]
 is finite.
  Given any Hilbert polynomial $L \in \mc{P}$, denote by 
  $S_L$ the constructible set obtained by taking the union 
  of the connected components $W$ of $T_i$ such that for 
  any $w \in W$, the Hilbert polynomial of $\det(\mc{F}_w)$
  equals $L$, for $1 \le i \le N$.
  
  By Theorem \ref{c2}, determinant is upper semi-continuous. One can easily to check that 
  for any given Hilbert polynomial $L_o \in \mc{P}$,
  the locus of points $s \in S$ such that $\det(\mc{F}_s)$
  has Hilbert polynomial strictly less than $L_o$, forms 
  an open subscheme. More precisely, the union of $S_{L'}$
  as $L'$ varies over all Hilbert polynomials strictly less 
  than $L_o$, is an open subscheme of $S$.
  In particular, for any $L \in \mc{P}$, $S_L$ is a closed 
  subscheme of the finite intersection of the open subschemes:
  \[\bigcap\limits_{L_i > L} \left( \coprod_{L'<L_i} S_{L'}
   \right).  \]
 Therefore, $S_L$ is in fact a locally closed subscheme of $S$.
  Hence, $\{S_L\}_{L \in \mc{P}}$ gives a stratification
  of $S$ satisfying conditions $(1)$ and $(2)$ of the theorem.
  
  We now prove $(3)$. Denote by $\mc{F}_{S'}$ the pull-back
  of $\mc{F}_S$ to $\mc{X}_{S'}$ under the morphism 
  $g_{\mc{X}}$. As before, define the set \[\mc{P}_{S'}:=
    \{L(\mc{F}_{s'})| s' \in S'\},      \]
where $L(\mc{F}_{s'})$ denotes the Hilbert polynomial of 
the determinant of $\mc{F}_{s'}$. Of course, $\mc{P}_{S'}
\subset \mc{P}$.

Since $g_{\mc{X}}^*\det_n(\mc{F}_S)$ is $(S_2)$-relative
to $\delta'$ and coincides with $\det_n(\mc{F}_{S'})$
over $U_{S'}:=U \times_S S'$, \cite[Proposition $3.6$]{hass}
implies that $g_{\mc{X}}^*\det_n(\mc{F}_S) \cong 
\det_n(\mc{F}_{S'})$. Indeed, the complement of $U_{S'}$ in 
$\mc{X}_{S'}$ is fiberwise of codimension at least $2$.
Furthermore, let $W$ be a connected component of $S'$.
By Lemma \ref{lem:detconst}, for any $w \in W$ 
the Hilbert polynomial, say $L_W$, of $\det(\mc{F}_w)$ 
remains fixed. Denote by $W'$ the fiber product $S_{L_W} \times_S W$
and $i: W' \to W$ the induced morphism. 
It is easy to check that the morphism $i: W' \hookrightarrow
W$ is bijective. Then, by \cite[Theorem $14.9$]{har}, the 
morphism $i$ is an isomorphism i.e., $W$ factors through 
$S_{L_W}$. Since $W$ is an arbitrary connected 
component of $S'$, this proves $(3)$, hence the theorem.
  \end{proof}

\begin{rem}
 Note that such a stratification $\{S_L\}_{L \in \mc{P}}$ 
 is unique (up to non-reduced scheme structure) in the sense 
 that if $\{S'_L\}_{L \in \mc{P}}$ is another stratification
 of $S$, then $(S_L)_{\red}=(S'_L)_{\red}$.
\end{rem}

 \subsection{Compactified Picard functor}\label{sec:pic}
 Suppose that $S$ is connected and there exists a section
 $\delta': S \to \mc{X}$ to the morphism $\delta$ such that 
 the image of $\delta'$ is contained in the smooth locus 
 of the morphism $\delta$.
 Fix a Hilbert polynomial $L$ of a 
 rank $1$, torsion-free sheaf on some fiber of $\delta$.
 Define the \emph{compactified Picard functor}, 
 denoted \[\ov{\mc{P}ic}_{\delta,L}: 
 (\mr{Sch}/S)^o \to (\mr{sets}),\]
 as the functor taking an $S$-scheme $T$ to the space of
 equivalent 
 classes of $T$-flat sheaves $\mc{L}_T$ on 
 $\mc{X} \times_S T$ with the following property: 
 for every geometric point $t \in T$, the 
 restriction $\mc{L}_t$ of $\mc{L}_T$ to $\mc{X}_t$ is a 
 torsion-free, rank $1$
 sheaf on $\mc{X}_t$ with Hilbert polynomial $L$. 
 The equivalence relation is given as follows: 
 two sheaves $\mc{L}_T$ and $\mc{L}'_T$ as above are 
 \emph{equivalent} if there exists an invertible sheaf $\mc{M}$ on $T$ such that 
 $\mc{L}_T \cong \mc{L}'_T \otimes p^*\mc{M}$ for the 
 natural projection morphism $p$ from $\mc{X} \times_S T$ to 
 $T$. By \cite[Theorem $3.3$]{lant1}, there exists
 a projective $S$-scheme 
 $\ov{\mr{Pic}}_{\delta,L}$, known as the \emph{relative compactified 
 Picard scheme}, that represents the 
 compactified Picard functor
 $\ov{\mc{P}ic}_{\delta,L}$.

Using Theorem \ref{th:str} we can now give a canonical determinant morphism
from the stratification of $S$ to a union of relative compactified Picard schemes.
 
 \begin{cor}\label{det11}
  There exists an 
  \emph{unique} morphism, which we call the \emph{determinant morphism}, 
  \[\mr{det}_S: \coprod_{L \in \mc{P}} S_L \to 
  \coprod_{L \in \mc{P}} \ov{\mr{Pic}}_{\delta,L}\]
  such that $\det_S$ takes the stratum $S_L$ to 
  $ \ov{\mr{Pic}}_{\delta,L}$ and 
    pull-back of the universal sheaf on 
  $\ov{\mr{Pic}}_{\delta,L}$
  to $S_L$ is isomorphic to $\det_n(\mc{F}_{S_L})$.
  Furthermore, for every $s \in S$,  $\mr{det}_S(s)$ is the 
  point corresponding to $\det(\mc{F}_s)$.
 \end{cor}

 \begin{proof}
  This is an immediate consequence of the 
  representability of $\ov{\mc{P}ic}_{\delta,L}$ 
  (see \cite[Theorem $3.3$]{lant1}) and Theorem \ref{th:str}.
  This proves the corollary.
 \end{proof}

 \section{Moduli spaces of sheaves with fixed determinant}
 In this section we construct a determinant morphism for the stratification of the moduli 
 space of stable sheaves over a normal, projective variety (Proposition \ref{cn31}). 
 This determinant morphism coincides with the classical one in 
 the case when the underlying variety is non-singular (in which case
 the stratification is trivial).
 In \S \ref{cn44}, we  
 define the moduli functor for stable 
 sheaves with fixed determinant. 
 Finally, we prove the 
 corepresentability of the moduli functor in certain cases (Theorem \ref{c13}). 
 We also show that if the moduli functor 
 is corepresentable then it arises 
 as the fiber of the determinant morphism.

 \subsection{Setup}\label{se:mod}
 Throughout this section $k$ is an algebraically closed field of arbitrary characteristic. Let $X$ be a normal, projective $k$-variety.
 Fix a Hilbert polynomial $L$ of a rank $1$ and a torsion-free sheaf $\mc{L}$ on $X$. Fix a Hilbert polynomial $P$ of a rank $n$, degree $d$ torsion-free sheaf on $X$ with $n$ coprime to $d$.

 \subsection{Moduli space of semi-stable sheaves}
   Recall, the moduli functor $\mc{M}_X(P)$ corresponding to semi-stable sheaves on $X$ with Hilbert polynomial $P$ is defined as follows:
\[\mc{M}_{X}(P) : (\mr{Sch}/k)^o \rightarrow (\mr{Sets})\] 
such that for a $k$-scheme $T$,         
   \[ \mc{M}_{X}(P)(T):= \left\{ \begin{array}{l}
    \mbox{ isomorphism classes of } \mbox{ coherent sheaves } \mc{F} \mbox{ on } X \times_k T,\, \mbox{ flat}\\
    \mbox{ over }T \mbox{ and for every point } {t} \mbox{ of } T, \, \mc{F}|_{X_t} \mbox{ is a torsion-free} \\ 
    \mbox{ semi-stable sheaf with Hilbert polynomial } P \mbox{ on } X_{t}
    \end{array} \right\}/ \sim \]
where $\mc{F} \sim \mc{G}$ if there exists an invertible sheaf $\mc{L}$ on $T$ such that $\mc{F} \cong \mc{G} \otimes p_T^*\mc{L}$
for the natural projection morphism  $p_T:X_T \to T$.

Since $X$ is integral, \cite[Lemma $1.2.13$ and $1.2.14$]{huy} implies that a semi-stable sheaf on $X$ of coprime rank and degree is also stable. 
Then  by \cite[Theorem $4.3.4$]{huy} the functor $\mc{M}_{X}(P)$ is universally corepresentable by a projective $k$-scheme $M_X(P)$.

\subsection{The construction of the moduli space}\label{sec:const}
We recall the construction of the moduli space $M_X(P)$ from \cite{huy}.
By \cite[Theorem $3.3.7$]{huy}, there exists an integer $e$ such that any semi-stable sheaf on $X$ with Hilbert polynomial $P$ is $e$-regular (in the sense of Castelnuovo-Mumford regularity). 
 Fix such an integer $e$. Denote by $\mc{H}:=\mo_{X}(-e)^{\oplus P(e)}$ and by $\mr{Quot}_{\mc{H}/X/P}$ the Quot scheme parameterizing all quotients of the form 
$\mc{H} \twoheadrightarrow \mc{Q}_0$, where $\mc{Q}_0$ has Hilbert polynomial $P$ (see \cite[\S $4.4$]{S1} for more details).

 Let $\mc{R}$ be the subset of $\mr{Quot}_{\mc{H}/X/P}$ consisting of all points corresponding to 
 quotients of the form $\mc{H}\twoheadrightarrow \mc{Q}_0$
 such that $\mc{Q}_0$ is semi-stable and 
 $H^0(\mc{Q}_0(e))$ is (non-canonically) isomorphic to $k^{\oplus P(e)}$.
   By \cite[Proposition $2.3.1$]{huy}, $\mc{R}$ is an open subscheme in $\mr{Quot}_{\mc{H}/X/P}$.   
   The group $\mr{GL}(P(e))=\mr{Aut}(\mc{H})$ acts on $\mr{Quot}_{\mc{H}/X/P}$ from the right by the composition $[\rho] \circ g=[\rho \circ g]$, where 
      $[\rho:\mc{H} \to \mc{F}] \in \mr{Quot}_{\mc{H}/X/P} \mbox{ and } g \in \mr{GL}(P(e))$. By \cite[Theorem $4.3.3$]{huy}, $\mc{R}$ is the set of semi-stable points of $\mr{Quot}_{\mc{H}/X/P}$
   under this group action. The moduli scheme 
   $M_X(P)$ of semi-stable sheaves on $X$ with Hilbert polynomial $P$ is the geometric quotient of $\mc{R}$ 
   under this action. Denote by 
   \begin{equation}\label{cn14}
    \pi:\mc{R} \to M_X(P)
   \end{equation}
 the corresponding quotient morphism. 
 By \cite[Corollary $4.3.5$]{huy}, the quotient $\pi$ is a $\mr{PGL}(P(e))$-bundle. 
 
 \begin{note}\label{det19}
  Denote by $\mc{Q}$ the universal
   quotient on $X \times \mr{Quot}_{\mc{H}/X/P}$ associated to $\mr{Quot}_{\mc{H}/X/P}$.
   For a locally closed subscheme $B \subset \mc{R}$, denote by $\mc{Q}_B$ the restriction $\mc{Q}|_{X \times B}$.
 \end{note}

 \begin{prop}\label{cn31}
  There exist a finite set of Hilbert polynomials 
$L_i$ of rank $1$, torsion-free sheaves and a stratification
$\coprod_i M_X(P,L_i)$ of $M_X(P)_{\red}$ 
such that we have a (canonical) 
determinant morphism 
\[\mr{det}_{M_X(P)_{\red}}: \coprod_i M_X(P,L_i)
\to \coprod_i \ov{\mr{Pic}}_{X,L_i}\]
mapping $M_X(P,L_i)$ to $\ov{\mr{Pic}}_{X,L_i}$. In particular, 
for any geometrically closed point $v \in M_X(P,L_i)$,
the determinant $\mr{det}(\mc{F}_v)$ of the 
associated semi-stable sheaf $\mc{F}_v$ has Hilbert polynomial $L_i$
and $\mr{det}_{M_X(P)_{\red}}$ maps $v$ to the point corresponding to $\det(\mc{F}_v)$.
 \end{prop}

\begin{proof}
Denote by $M_X(P)_{\red}$ (resp. $\mc{R}_{\red}$) 
the reduced scheme associated to $M_X(P)$ (resp. $\mc{R}$).
Since $\pi$ is a smooth morphism,
we get the following fiber product 
diagram:
\[\begin{diagram}
   \mc{R}_{\red}&\rTo&\mc{R}\\
   \dTo^{\widetilde{\pi}}&\square&\dTo^{\pi}\\
   M_X(P)_{\red}&\rTo&M_X(P)
     \end{diagram}\]
where the two horizontal maps are natural.
Since $\pi$ is an universal geometric quotient, we have 
$\widetilde{\pi}$ is also an universal geometric quotient.
By Theorem \ref{th:str} and Corollary \ref{det11}, 
there exist a finite set of Hilbert polynomials 
$L_i$ of rank $1$, torsion-free sheaves and a stratification
$\coprod_i \mc{R}_i$ of $\mc{R}_{\red}$ 
such that we have an unique 
determinant morphism 
\[\mr{det}_{\mc{R}_{\red}}: \coprod_i \mc{R}_i
\to \coprod_i \ov{\mr{Pic}}_{X,L_i}\]
such that $\mc{R}_i$ maps to $\ov{\mr{Pic}}_{X,L_i}$ and
 the pull-back of the universal sheaf on 
  $\ov{\mr{Pic}}_{X,L_i}$
  to $\mc{R}_i$ is isomorphic to 
  $\det_n(\mc{Q}_{\mc{R}_i})$. Moreover for every $y \in \mc{R}$, 
  $\mr{det}_{\mc{R}_{\red}}(y)$ is the 
  point corresponding to $\det(\mc{Q}_y)$.
  Furthermore, the $\mr{GL}(P(e))$-action on $\mc{R}$
  restricts to a $\mr{GL}(P(e))$ action on $\mc{R}_i$, for 
  each $i$. Taking the quotient by the $\mr{GL}(P(e))$-action,
  we obtain a locally closed subscheme $M_X(P,L_i)$ of 
  $M_X(P)$ such that the restriction $\pi_{L_i}$ of $\pi$
  to $\mc{R}_i$ is the universal geometric quotient 
  \[\pi_{L_i}: \mc{R}_i \to M_X(P,L_i).\]
  In fact, $\pi_i$ is a $\mr{PGL}(P(e))$-bundle.
 By the universal property of corepresentable functors,
 $\det_{\mc{R}_{\red}}$ gives rise to a morphism 
 \[\mr{det}_{M_X(P)_{\red}}: \coprod_i {M}_X(P,L_i) \to 
  \coprod_i \ov{\mr{Pic}}_{X,L_i} \]
  such that $\det_{M_X(P)_{\red}} \circ \widetilde{\pi}=
  \det_{\mc{R}_{\red}}$.
  This proves the proposition.
  \end{proof}
  
  \begin{note}  
 Let $o$ be a closed point in $\coprod \ov{\mr{Pic}}_{X,L_i}$.
 Denot by $\mc{L}$ the rank $1$, torsion-free sheaf on $X$ corresponding
 to the point $o$. Denote by 
 \[M_X(P,\mc{L}):=\mr{det}_{M_X(P)_{\red}}^{-1}(o)\, 
  \mbox{ and } \mc{R}_{\mc{L}}:=
  \mr{det}_{\mc{R}_{\red}}^{-1}(o). \]
Note that, for any $i$, the pre-image of 
$\ov{\mr{Pic}}_{X,L_i}$
under the morphism $\det_{M_{X}(P)_{\red}}$ is
the stratum $M_X(P,L_i)$ 
of $M_X(P)_{\red}$. As a result, 
$M_X(P,\mc{L})$ is a locally 
closed subscheme of $M_X(P)$. 
Hence, it is a quasi-projective $k$-scheme. 
Furthermore, the morphism 
$\pi$ restricts to the morphism \[\pi_{\mc{L}}: \mc{R}_{\mc{L}} 
    \to M_X(P,\mc{L})                             \]
which makes $\mc{R}_{\mc{L}}$ a $\mr{PGL}(P(e))$-bundle over 
$M_X(P,\mc{L})$.
  \end{note}

  \subsection{Moduli functor with fixed determinant}\label{cn44}
  Let $\mc{L}$ be a rank $1$ semi-stable sheaf on $X$. Define the \emph{moduli functor of semi-stable sheaves with fixed alternating determinant}:
     \[\mc{M}_{X}(P,\mc{L}) : (\mr{Sch}/k)^o \rightarrow (\mr{Sets})\] 
such that for a $k$-scheme $T$,         
   \[ \mc{M}_{X}(P,\mc{L})(T):= \left\{ \begin{array}{l}
    \mbox{ isomorphism class of sheaves } \mc{F}_T \in \mc{M}_X(P)(T) \mbox{ such that }\\
     \det_n(\mc{F}_T) \cong \mc{L} \otimes_k \mo_T \mbox{ and for all } t \in T,  \, \det(\mc{F}_t) \cong \mc{L} \otimes_k k(t)
    \end{array} \right\}/ \sim \] 
 where $\mc{F}_T \sim \mc{G}_T$ if there exists an invertible sheaf $\mc{N}$ on $T$ such that 
 $\mc{F}_T \cong \mc{G}_T \otimes p_T^*\mc{N}$ for the natural projection $p_T:X_T \to T$. By Lemma \ref{lem:func}, we have $\mc{M}_X(P,\mc{L})$
 is in fact a functor.
 
 The proof of the following result is a modification of the 
 proof of the classical corepresentability result given in 
 \cite[Lemma $4.3.1$]{huy}.
 
 \begin{thm}\label{c13}
 We have the following:
 \begin{enumerate}
  \item Suppose $\mc{M}_X(P,\mc{L})$ is corepresentable by a locally closed subscheme of $M_X(P)$, which we denote by $N_{\mc{L}}$. Then $N_{\mc{L},\red}=M_X(P,\mc{L})_{\red}$.
 \item Assume that $M_X(P)$ is reduced and $L$ is the smallest 
 Hilbert polynomial of $\det(\mc{Q}_s)$ as $s$ 
 varies over $\mc{R}$.
 Then for any rank $1$, $(S_2)$-sheaf $\mc{L}$ on $X$
 with Hilbert polynomial $L$, we have 
 $\mc{M}_X(P,\mc{L})$ is 
 universally corepresentable by the quasi-projective 
 $k$-scheme  $M_X(P,\mc{L})$.
 \end{enumerate}
   \end{thm}

 \begin{proof}
Let $S$ be a $k$-scheme and $\mc{F}_S \in \mc{M}_X(P,\mc{L})(S)$. 
Following notations as in \S \ref{sec:const}, $\mc{F}_s$ is 
$e$-regular for every $s \in S$. Denote by $p$ (resp. $q$)
the natural morphism from $X \times_k S$ to $S$ (resp. $X$).
By \cite[Proposition $2.1.2$]{huy}, 
\[\mc{G}:=p_*(\mc{F}_S \otimes q^*\mo_X(e))\]
is locally-free of rank $P(e)$. There exists a canonical surjective 
morphism from $p^*\mc{G}$ to $\mc{F}_S \otimes q^*\mo_X(e)$ 
(see proof of \cite[Theorem III.$8.8$]{R1}). Tensoring this morphism
by $q^*\mo_X(-e)$, we get the surjective morphism:
\[\phi: p^*\mc{G} \otimes q^*\mo_X(-e) \to \mc{F}_S.\]
 Denote by $R(\mc{F}_S):=\mr{Isom}(k^{\oplus P(e)} \otimes_k \mo_S,
 \mc{G})$ the frame bundle associated to $\mc{G}$ (see 
 \cite[$4.2.3$]{huy}) parameterizing fiber-wise isomorphisms from 
 $k^{\oplus P(e)} \otimes_k \mo_S$ to $\mc{G}$.
 Let $\pi_S: R(\mc{F}_S) \to S$ be the natural projection and 
 $\nu: k^{\oplus P(e)} \otimes_k \mo_S \otimes_{\mo_S} \mo_{R(\mc{F}_S)}
 \to \mc{G} \otimes_{\mo_S} \mo_{R(\mc{F}_S)}$ the universal 
 trivialization of $\mc{G}$ over $R(\mc{F}_S)$. 
 Pulling back the morphism $\phi$ to $X \times R(\mc{F}_S)$ via 
 the morphism $\mr{id} \times \pi_S$ and composing with the morphism 
 $\nu$, we get the following surjective morphism:
 \[\widetilde{q}_{_{\mc{F}_S}}:\mc{H} \otimes_k \mo_{R(\mc{F}_S)} \twoheadrightarrow \mc{F}_S \otimes_{\mo_S} \mo_{R(\mc{F}_S)}.\]
 By the universal property of the Quot scheme, this gives rise to a morphism 
\[\widetilde{\Phi}_{\mc{F}_S}:R(\mc{F}_S) \to \mr{Quot}_{\mc{H}/X/P}\]
which factors through $\mc{R}$. 
If $u$ is a geomtric point such that the corresponding 
semi-stable sheaf $\mc{F}_u$ over $X_u$ has determinant different
from $\mc{L}$, then $\widetilde{\Phi}_{\mc{F}_S}$ factors through 
the complement $\mc{R}_u^c:=\mc{R}\backslash \pi^{-1}(u)$.
Note that by the coprime assumption of degree and rank 
that semi-stability coincides with stability. Suppose $\mc{R}$ is reduced and $L$ is the smallest
Hilbert polynomial of the determinant of $\mc{Q}_s$ as 
$s$ varies over $\mc{R}$. Then, by
Corollary \ref{det11}  $\widetilde{\Phi}_{\mc{F}_S}$ factors
through $\mc{R}_{\mc{L}}$.

Note that, $R(\mc{F}_S)$ is a $\mr{GL}(P(e))$-bundle over $S$ and 
$\widetilde{\Phi}_{\mc{F}_S}$ is $\mr{GL}(P(e))$-equivariant.
Hence the morphism $\widetilde{\Phi}_{\mc{F}_S}$ induces a natural transformation 
\[\Phi_{\mc{F}_S}^c:\underline{R(\mc{F}_S)}/\underline{\mr{GL}(P(e))} \to \underline{\mc{R}^c_u}/\underline{\mr{GL}(P(e))}.\]
Note that $\mc{R}_{\mc{L}}$ is equipped with a natural 
 $\mr{GL}(P(e))$-action on the left, acting on the free $\mo_X$-module
 $\mc{H}$, as described in \S \ref{sec:const}. 
Therefore, if $\mc{R}$ is reduced, $\widetilde{\Phi}_{\mc{F}_S}$ induces a natural transformation 
\[\Phi_{\mc{F}_S}:\underline{R(\mc{F}_S)}/\underline{\mr{GL}(P(e))} \to \underline{\mc{R}_{\mc{L}}}/\underline{\mr{GL}(P(e))}.\]
 Since $\pi_S$ is a $\mr{GL}(P(e))$-bundle morphism, 
  the functor $\underline{R(\mc{F}_S)}/\underline{\mr{GL}(P(e))}$ is corepresented by $\underline{S}$.
 Using the universal property of corepresentable functors, the 
 natural transformations $\Phi_{\mc{F}_S}^c$ (resp. $\Phi_{\mc{F}_S}$)
 induce a natural transformation from $\underline{S}$ to 
 $\underline{\mc{R}^c_u}/\underline{\mr{GL}(P(e))}$ (resp. 
 $\underline{\mc{R}_{\mc{L}}}/\underline{\mr{GL}(P(e))}$).
 Then, the identity map
 $\mr{id}_S \in \underline{S}(S)$ defines an element in 
 $\underline{\mc{R}^c_u}/\underline{\mr{GL}(P(e))}(S)$
 (resp. $\underline{\mc{R}_{\mc{L}}}/\underline{\mr{GL}(P(e))}(S)$).
 This defines a natural transformation:
 \[\tau_{0}:\mc{M}_X(P,\mc{L}) \to \underline{\mc{R}^c_u}/\underline{\mr{GL}(P(e))}\]
 and in the case $\mc{R}$ is reduced and $L$ is minimal
 as mentioned above, we get a natural transformation
 \[\tau_{0}^{\mr{red}}:\mc{M}_X(P,\mc{L}) \to 
  \underline{\mc{R}_{\mc{L}}}/\underline{\mr{GL}(P(e))}. \]
The universal family over $\mc{R}_{\mc{L}}$ gives the inverse natural 
transformation from $\underline{\mc{R}_{\mc{L}}}/\underline{\mr{GL}(P(e))}$
to $\mc{M}_X(P,\mc{L})$. Since 
$\underline{\mc{R}_{\mc{L}}}/\underline{\mr{GL}(P(e))}$ is corepresented
by $M_X(P,\mc{L})$, we conclude that in the case $\mc{R}$ is reduced and $L$ is minimal,
$\mc{M}_X(P,\mc{L})$ is corepresented by $M_X(P,\mc{L})$. This proves
$(2)$.

Suppose (without assumption on reducedness on $\mc{R}$) that 
$\mc{M}_X(P,\mc{L})$ is corepresentable by a locally closed subscheme of $M_X(P)$,
which we call $\mc{N}_{\mc{L}}$ .
   Note that $M_X(P) \backslash \{u\}$ corepresents the categorical quotient $\underline{\mc{R}^c_u}/\underline{\mr{GL}(P(e))}$.
 Since $N_{\mc{L}}$ corepresents $\mc{M}_X(P,\mc{L})$, the universal property of corepresentable functors implies that there exists a natural transformation $\ov{\tau}_0$ from 
 $\underline{\mc{N}_{\mc{L}}}$ to $ \underline{M_X(P)\backslash \{u\}}$
 making the following diagram commute:
 \[\begin{diagram}
  \mc{M}_X(P,\mc{L})&\rTo^{\tau_0}&\underline{\mc{R}^c_u}/\underline{\mr{GL}(P(e))}&\rTo^{\psi}&\underline{\mc{R}}/\underline{\mr{GL}(P(e))}\\
  \dTo&\circlearrowleft&\dTo&\circlearrowleft&\dTo\\
  \underline{N_{\mc{L}}}&\rTo^{\ov{\tau_0}}&\underline{M_X(P)\backslash \{u\}}&\rTo^{\ov{\psi}}&\underline{M_X(P)}
   \end{diagram}\]
   where the morphisms $\psi$ and $\ov{\psi}$ are natural.
 By Yoneda lemma, this gives rise to a morphism $\ov{\psi} \circ \ov{\tau_0}:N_{\mc{L}} \to M_X(P)$. Since $N_{\mc{L}}$ is a locally closed subscheme of $M_X(P)$,
 $\ov{\psi} \circ \ov{\tau_0}$ is a locally closed immersion (use the uniqueness part of the universal property of categorical quotient). Hence, 
 $\ov{\tau_0}:N_{\mc{L}} \to M_X(P) \backslash \{u\}$ is a locally closed immersion i.e., $u \not\in N_{\mc{L}}$. In other words, 
 for any geometric point $u \in M_X(P)\backslash M_X(P,\mc{L})$,
 we have $u \not\in \mc{N}_{\mc{L}}$.
 Denote by $W$ the (scheme-theoretic) 
 intersection of $M_X(P,\mc{L})$ and $N_{\mc{L}}$, 
 as locally closed subschemes in $M_X(P)$. Let 
 $\phi:W \hookrightarrow  N_{\mc{L}}$ denote the natural inclusion.
 The above argument implies that 
 the morphism $\phi$ is 
 surjective, thereby also bijective.
 Given a geometric point $u \in M_X(P)$, $u \in M_X(P,\mc{L})$ if and only if for the corresponding isomorphism class of 
a semi-stable sheaf $\mc{F}_u$ on $X$, $\det(\mc{F}_u) \cong \mc{L}$.
 Since $N_{\mc{L}}$ corepresents $\mc{M}_X(P,\mc{L})$, all such geometric points $u$ must also be contained 
 in $N_{\mc{L}}$. Hence, $\phi_1: W \hookrightarrow M_X(P,\mc{L})$
 is bijective. Using \cite[Theorem $14.9$]{har} we conclude by the  
  bijectivity of $\phi$ and $\phi_1$  that 
  \[N_{\mc{L},\red} \cong W_{\red} \cong M_X(P,\mc{L})_{\red}.\]
  This proves $(1)$ and hence the theorem.
 \end{proof}

\section{Applications to geometry}
In this section we show how the stratification given in 
Theorem \ref{th:str} can be used to study certain kinds of singularities of 
curves with fixed Hilbert polynomial (Theorem \ref{th:geom}).
As an example, we apply this stratification
to study the   Hilbert scheme of curves in $\p3$ with Hilbert polynomial $3n+1$ originally studied by Schlessinger and Piene
(Example \ref{ex:sch}).

\subsection{Flag Hilbert schemes}\label{sec:flag}
Let $P$ (resp. $Q$) be a Hilbert polynomial of degree $1$
(resp. $2$). Denote by $\mr{Hilb}_{P,Q}$ the flag 
Hilbert scheme parameterizing pairs $(C \subset X)$,
where $C$ (resp. $X$) is closed subscheme of $\p3$ with 
Hilbert polynomial $P$ (resp. $Q$). See 
\cite[\S $4.5$]{S1} for a detailed study. We are 
interested in the sub-locus $\mr{Hilb}^{\mr{nor}}_{P,Q}
\subset \mr{Hilb}_{P,Q}$ parameterizing pairs 
$(C \subset X)$ with $X$ a normal hypersurface in $\p3$.
By \cite[Theorem $12.2.4$]{ega43}, we observe that 
$\mr{Hilb}^{\mr{nor}}_{P,Q}$ is an open subscheme in 
$\mr{Hilb}_{P,Q}$. Let us consider the universal 
family restricted to  
$\mr{Hilb}_{P,Q}^{\mr{nor}}$ i.e., flat, projective 
morphisms \[\pi_1: \mc{C} \to \mr{Hilb}_{P,Q}^{\mr{nor}},\,
\pi_2: \mc{X} \to  \mr{Hilb}_{P,Q}^{\mr{nor}}\]
 and closed immersions $i:\mc{C} \hookrightarrow \mc{X}$
 and $j: \mc{X} \hookrightarrow \p3 \times \mr{Hilb}_{P,Q}^{\mr{nor}}$
such that $\pi_2 \circ i = \pi_1$ and $\pr_2 \circ j=\pi_2$,
where $\pr_2$ is the natural morphism from $\p3 \times \mr{Hilb}_{P,Q}^{\mr{nor}}$
to $\mr{Hilb}_{P,Q}^{\mr{nor}}$.
The closed immersion $i$ then 
gives rise to a surjective morphism from $\mo_{\mc{X}}$ 
to $\mo_{\mc{C}}$. Denote by $\I_{\mc{C}}$ the kernel 
of this morphism i.e., we have a short exact sequence:
\begin{equation}\label{eqn17}
 0 \to \I_{\mc{C}} \to \mo_{\mc{X}} \to \mo_{\mc{C}} \to 0.
\end{equation}
For every $y \in \mr{Hilb}_{P,Q}^{\mr{nor}}$, denote by $\mc{C}_y:=
\pi_1^{-1}(y)$ and $\mc{X}_y:=\pi_2^{-1}(y)$. 
Note that, as $\mo_{\mc{C}}$ and $\mo_{\mc{X}}$ are 
$\mr{Hilb}_{P,Q}^{\mr{nor}}$-flat, so is $\I_{\mc{C}}$ (use 
the short exact sequence \eqref{eqn17}).
This implies for every 
$y \in \mr{Hilb}_{P,Q}^{\mr{nor}}$, we have 
$\I_{\mc{C}_y}$ is a sub-sheaf of $\mo_{\mc{X}_y}$. Hence it is torsion-free.

\subsection{Stratification of the flag Hilbert scheme}
Notations as \S \ref{sec:flag} above. For simplicity
of notation, denote by $S:=
(\mr{Hilb}_{P,Q}^{\mr{nor}})_{\red}$, the reduced 
scheme associated to $\mr{Hilb}_{P,Q}^{\mr{nor}}$.
Let \[\{S_L\}_{L \in \mc{P}}\] be the stratification
of $S$ stated in Theorem \ref{th:str} after replacing 
$\mc{F}_S$ by $\I_{\mc{C}}$, parameterized by the set 
$\mc{P}$ of Hilbert polynomials:
\[\mc{P}:=\{L(\I_{\mc{C}_s})| s \in S\},\]
where $L(\I_{\mc{C}_s})$ is the Hilbert polynomial of 
$\det(\I_{\mc{C}_s})$.

Note that, for a given $s \in S$, the corresponding curve
$\mc{C}_s$ need not be Cohen-Macaulay. We call the 
\emph{Cohen-Macaulayfication of} $\mc{C}_s$, denoted
$\mo_{\mc{C}_s}^{\mr{CM}}$, the $\mo_{\mc{X}_s}$-module:
\[\mo_{\mc{C}_s}^{\mr{CM}}:=\mc{E}xt^1_{\mc{X}_s}
(\mc{E}xt^1_{\mc{X}_s}(\mo_{\mc{C},s},\mo_{\mc{X},s}),
\mo_{\mc{X},s}).\]
Using \cite[Theorem $3.3.10$]{brun1}, one can check that 
if $\mc{C}_s$ is Cohen-Macaulay (for example, a Cartier
divisor in $\mc{X}_s$), then 
$\mo_{\mc{C}_s} \cong \mo_{\mc{C}_s}^{\mr{CM}}$.
Denote by $P(\mo_{\mc{C}_s}^{\mr{CM}})$ the 
Hilbert polynomial of $\mo_{\mc{C}_s}^{\mr{CM}}$. Let $\mc{M}$ denote the set $\{P(\mo_{\mc{C}_s}^{\mr{CM}})| s \in S\}$. We show:

\begin{thm}\label{th:geom}
 The set $\mc{M}$ is finite and there is a bijective map 
 (of sets) from $\mc{P}$ to $\mc{M}$, sending a Hilbert 
 polynomial $L \in \mc{P}$ to $Q-L \in \mc{M}$. 
 In particular,
\begin{enumerate}
  \item for every $L \in \mc{P}$, we have $L \ge Q-P$
  with equality if and only if there exists $s \in S_L$ 
  for which $\mc{C}_s$ is Cohen-Macaulay,
  \item $s \in S_L$ if and only if the 
  Cohen-Macaulayfication $\mo_{\mc{C}_s}^{\mr{CM}}$ of 
  $\mc{C}_s$ has Hilbert polynomial $Q-L$.
 \end{enumerate}
 \end{thm}

\begin{proof}
 For any $s \in S$, the short exact sequence \eqref{eqn17}
 restricts to 
 \begin{equation}\label{eqn20}
  0 \to \I_{\mc{C}_s} \to \mo_{\mc{X}_s} \to \mo_{\mc{C}_s}
  \to 0.
 \end{equation}
Dualizing this short exact sequence, we get 
\begin{equation}\label{eqn19}
 0 \to \mo_{\mc{X}_s} \to (\I_{\mc{C}_s})^{\vee} \to 
 \mc{E}xt^1_{\mc{X}_s}(\mo_{\mc{C}_s}, \mo_{\mc{X}_s})
 \to 0.
\end{equation}
Note that, as $(\I_{\mc{C}_s})^{\vee}$ is reflexive 
(see \cite[Corollary $1.2$]{stabhar}) and $\mc{X}_s$ is 
a surface, we have $(\I_{\mc{C}_s})^{\vee}$ is maximal 
Cohen-Macaulay. Hence, by \cite[Theorem $3.3.10$]{brun1},
we have $\mc{E}xt^1_{\mc{X}_s}
((\I_{\mc{C}_s})^\vee,\mo_{\mc{X}_s})=0$. Dualizing 
\eqref{eqn19} then gives us the short exact sequence:
\begin{equation}\label{eqn18}
 0 \to \det(\I_{\mc{C}_s}) \to \mo_{\mc{X}_s} \to 
 \mo_{\mc{C}_s}^{\mr{CM}} \to 0.
\end{equation}
Therefore,  
the Hilbert polynomial $P(\mo_{\mc{C}_s}^{\mr{CM}})$ 
of $\mo_{\mc{C}_s}^{\mr{CM}}$ equals 
$Q-L(\I_{\mc{C}_s})$.
Then, the finiteness of $\mc{P}$ proved in Theorem
\ref{th:str} immediately implies the finiteness of $\mc{M}$
and we have the bijective map from $\mc{P}$ to $\mc{M}$
as mentioned in the statement of the theorem.

Since $\I_{\mc{C}_s}$ is a torsion-free 
$\mo_{\mc{X}_s}$-module, there is a natural inclusion 
\begin{equation}\label{eqninc}
 \I_{\mc{C}_s} \hookrightarrow \det(\I_{\mc{C}_s})
\end{equation}
which is an isomorphism over a complement of a subvariety 
in $\mc{X}_s$ of codimension at least $2$.
Note that, the inclusion is an isomorphism if and only if 
$\I_{\mc{C}_s}$ is reflexive. Moreover, if $\I_{\mc{C}_s}$ is reflexive, then 
by depth comparison applied to \eqref{eqn20} we conclude that 
$\mc{C}_s$ is Cohen-Macaulay . 
Conversely, if $\mc{C}_s$
is Cohen-Macaulay, then once again by depth comparison
applied to \eqref{eqn20}, we have $\I_{\mc{C}_s}$ is 
reflexive. Hence, the inclusion in \eqref{eqninc} 
 is an isomorphism if and only 
if $\mc{C}_s$ is Cohen-Macaulay.
This implies 
that the Hilbert polynomial of $\I_{\mc{C}_s}$ is 
less than or equal to $\det(\I_{\mc{C}_s})$, with 
equality if and only if $\mc{C}_s$ is Cohen-Macaulay. 
By \eqref{eqn20},
the Hilbert polynomial of $\I_{\mc{C}_s}$ is $Q-P$.
 Therefore, for every $s \in S$, the 
 Hilbert polynomial of $\det(\I_{\mc{C}_s})$
is greater than or equal to $Q-P$ with equality 
 if and only if
 $\mc{C}_s$ is Cohen-Macaulay. 
 Moreover, using \eqref{eqn18}, we observe that 
 $s \in S_L$ if and only if the 
  Cohen-Macaulayfication $\mo_{\mc{C}_s}^{\mr{CM}}$ of 
  $\mc{C}_s$ has Hilbert polynomial $Q-L$.
 This proves the theorem.
\end{proof}

\subsection{Example}\label{ex:sch}
We now apply Theorem \ref{th:geom} to the classical 
example of Schlessinger and Piene \cite{piesch}.
Denote by $P(n)=3n+1$ and $Q_d$ the Hilbert polynomial
of a degree $d$ hypersurface in $\p3$. Schlessinger and Piene showed that the Hilbert scheme 
corresponding to the Hilbert polynomial $3n+1$ consists 
of two irreducible components, one
parameterizing smooth twisted cubic curves and the 
other component
parameterizing plane cubic curves along 
with an embedded point 
(i.e., the curve contains an embedded point and the 
underlying reduced curve is a plane cubic curve). Note that,
given a curve $C$ containing an embedded point, say at $o$,
the depth of $\mo_{C,o}$ is zero. Hence $C$ 
cannot be Cohen-Macaulay. For $d \gg 0$, by a Bertini type theorem (see \cite[Theorem $7$]{kleim}), 
there exist normal, projective hypersurfaces in $\p3$ 
containing a plane cubic curve along 
with an embedded point. For such a $d$, consider the 
Hilbert scheme $\mr{Hilb}_{P,Q_d}$ and the subscheme
$S:=\mr{Hilb}_{P,Q_d}^{\mr{nor}}$ where $X$ is normal.  The subscheme $S$ parameterizes
pairs $(C \subset X)$ with $C$ (resp. $X$) closed subschemes
in $\p3$ with Hilbert polynomial $P$ (resp. $Q_d$). 
Then, the set $\mc{P}$ as in Theorem \ref{th:str} consists of two 
quadratic polynomials:
\[P'(n):=Q_d(n)-3n \mbox{ and } P(n):=Q_d(n)-(3n+1)\]
and the corresponding stratum $S_{P}$ (resp. 
$S_{P'}$) parametrizes pairs $(C,X)$ with $C$ a twisted
cubic curve (resp. a cubic curve with an embedded point)
and $X$ is normal. Indeed, if $C$ is a plane cubic curve 
with an embedded point, we have the short
exact sequence:
\[0 \to K \to \mo_C \to \mo_{C_{\red}} \to 0\]
where $C_{\red}$ is the reduced plane cubic curve and $K$
is supported at the embedded point of $C$.
Suppose that $C$ is contained in a normal surface $X$ in 
$\p3$. Applying $\Hc_X(-,\mo_X)$ to the short exact 
sequence and using 
$\Hc_X(K,\mo_X)=0=\mc{E}xt^1_X(K,\mo_X)$ (see 
\cite[Theorem $3.3.10$]{brun1}) we get 
\[\mc{E}xt_X^1(\mo_{C_{\red}}, \mo_X) 
          \xrightarrow{\sim} \mc{E}xt^1_X(\mo_C,\mo_X). \]
Dualizing once again we get $\mo_{C_{\red}}^{\mr{CM}}
\cong \mo_C^{\mr{CM}}$. Since $C_{\red}$ is reduced, it 
is Cohen-Macaulay. Hence $\mo_{C_{\red}}^{\mr{CM}}$
is isomorphic to $\mo_{C_{\red}}$.
Note that, the Hilbert polynomial of $\mo_{C_{\red}}$ is 
$3n$. Therefore, the Hilbert polynomial of 
$\mo_C^{\mr{CM}}$ is $Q_d - P'$.


\begin{thebibliography}{10}

\bibitem{barthmod}
{\sc W.~Barth}, {\em Moduli of vector bundles on the projective plane},
  Inventiones mathematicae, 42 (1977), pp.~63--91.

\bibitem{genpre}
{\sc S.~Basu, A.~Dan, and I.~Kaur}, {\em Generators of the cohomology ring,
  after {N}ewstead}, Proceedings of the American Mathematical Society, 150
  (2022), pp.~2569--2577.

\bibitem{bayemmp}
{\sc A.~Bayer and E.~Macr{\`\i}}, {\em {MMP} for moduli of sheaves on {K3}s via
  wall-crossing: nef and movable cones, {L}agrangian fibrations}, Inventiones
  mathematicae, 198 (2014), pp.~505--590.

\bibitem{bhos2}
{\sc U.~N. Bhosle}, {\em Generalised parabolic bundles and applications to
  torsionfree sheaves on nodal curves}, Arkiv f{\"o}r Matematik, 30 (1992),
  pp.~187--215.

\bibitem{brun1}
{\sc W.~Bruns and H.~J. Herzog}, {\em Cohen-{M}acaulay rings}, Cambridge
  University Press, 1998.

\bibitem{casablow}
{\sc C.~Casagrande, G.~Codogni, and A.~Fanelli}, {\em The blow-up of
  {$\mathbb{P}^4$} at 8 points and its {F}ano model, via vector bundles on a
  del {P}ezzo surface}, Revista Matem{\'a}tica Complutense, 32 (2019),
  pp.~475--529.



\bibitem{mumf}
{\sc A.~Dan and I.~Kaur}, {\em Generalization of a conjecture of {M}umford},
  Advances in Mathematics, 383 (2021), p.~107676.

\bibitem{hodc}
\leavevmode\vrule height 2pt depth -1.6pt width 23pt, {\em Hodge conjecture for
  the moduli space of semi-stable sheaves over a nodal curve}, Annali di
  Matematica Pura ed Applicata (1923 -), DOI- 10.1007/s10231-022-01225-7
  (2022).

\bibitem{donpoly}
{\sc S.~K. Donaldson}, {\em Polynomial invariants for smooth four-manifolds},
  Topology, 29 (1990), pp.~257--315.

\bibitem{elli3}
{\sc G.~Ellingsrud, R.~Piene, and S.~A. Str{\o}mme}, {\em On the variety of
  nets of quadrics defining twisted cubics}, in Space curves, Springer, 1981,
  pp.~84--96.

\bibitem{elli4}
{\sc G.~Ellingsrud and S.~Str{\"{o}}mme}, {\em Bott's formula and enumerative
  geometry}, Journal of the American Mathematical Society, 9 (1996),
  pp.~175--193.

\bibitem{elstr}
{\sc G.~Ellingsrud and S.~A. Str{\o}mme}, {\em On the rationality of the moduli
  space for stable rank-2 vector bunoles on {$\mathbb{P}^2$}}, in
  Singularities, representation of algebras, and vector bundles, Springer,
  1987, pp.~363--371.

\bibitem{giesli}
{\sc D.~Gieseker and J.~Li}, {\em Irreducibility of moduli of rank-2 vector
  bundles on algebraic surfaces}, Journal of Differential Geometry, 40 (1994),
  pp.~23--104.

\bibitem{giesli2}
\leavevmode\vrule height 2pt depth -1.6pt width 23pt, {\em Moduli of high rank
  vector bundles over surfaces}, Journal of the American Mathematical Society,
  9 (1996), pp.~107--151.

\bibitem{gortz}
{\sc U.~G{\"o}rtz and T.~Wedhorn}, {\em Algebraic Geometry}, Springer, 2010.

\bibitem{ega43}
{\sc A.~Grothendieck}, {\em {\'E}l{\'e}ments de g{\'e}om{\'e}trie
  alg{\'e}brique (r{\'e}dig{\'e}s avec la collaboration de {J}ean
  {D}ieudonn{\'e}): {IV}. {\'e}tude locale des sch{\'e}mas et des morphismes de
  sch{\'e}mas, troisi{\`e}me partie}, Publications Math{\'e}matiques de
  l'IH{\'E}S, 28 (1966), pp.~5--255.

\bibitem{har}
{\sc J.~Harris}, {\em Algebraic geometry: a first course}, vol.~133, Springer
  Science \& Business Media, 2013.

\bibitem{R1}
{\sc R.~Hartshorne}, {\em Algebraic Geometry}, Graduate text in Mathematics-52,
  Springer-Verlag, 1977.

\bibitem{stabhar}
\leavevmode\vrule height 2pt depth -1.6pt width 23pt, {\em Stable reflexive
  sheaves}, Mathematische Annalen, 254 (1980), pp.~121--176.

\bibitem{hass}
{\sc B.~Hassett and S.~J. Kov{\'a}cs}, {\em Reflexive pull-backs and base
  extension}, Journal of Algebraic Geometry, 13 (2004), pp.~233--248.

\bibitem{hulek1}
{\sc K.~Hulek}, {\em Stable rank-2 vector bundles on {$\mathbb{P}^2$} with
  {$c_1$} odd}, Mathematische Annalen, 242 (1979), pp.~241--266.

\bibitem{huy}
{\sc D.~Huybrechts and M.~Lehn}, {\em The geometry of moduli spaces of
  sheaves}, Springer, 2010.

\bibitem{ind}
{\sc I.~Kaur}, {\em The ${C_1}$ conjecture for the moduli space of stable
  vector bundles with fixed determinant on a smooth projective curve,}.
\newblock \url{https://refubium.fu-berlin.de/handle/fub188/6611}, 2016.
\newblock Ph.D. thesis, Freie University Berlin.

\bibitem{K4}
\leavevmode\vrule height 2pt depth -1.6pt width 23pt, {\em Existence of
  semistable vector bundles with fixed determinants}, Journal of Geometry and
  Physics, 138 (2019), pp.~90--102.

\bibitem{ind2}
\leavevmode\vrule height 2pt depth -1.6pt width 23pt, {\em A pathological case
  of the {$C_1$} conjecture in mixed characteristic}, Mathematical Proceedings
  of the Cambridge Philosophical Society, 167 (2019), pp.~61--64.

\bibitem{king}
{\sc A.~King and A.~Schofield}, {\em Rationality of moduli of vector bundles on
  curves}, Indagationes Mathematicae, 10 (1999), pp.~519--535.

\bibitem{kleim}
{\sc S.~L. Kleiman and A.~B. Altman}, {\em Bertini's theorem for hypersurface
  sections containing a subscheme}, Communications in Algebra, 7 (1979),
  pp.~775--790.

\bibitem{langmix}
{\sc A.~Langer}, {\em Moduli spaces of sheaves in mixed characteristic}, Duke
  Mathematical Journal, 124 (2004), pp.~571--586.

\bibitem{langa}
\leavevmode\vrule height 2pt depth -1.6pt width 23pt, {\em Semistable sheaves
  in positive characteristic}, Annals of mathematics,  (2004), pp.~251--276.

\bibitem{langc}
\leavevmode\vrule height 2pt depth -1.6pt width 23pt, {\em Moduli spaces and
  {C}astelnuovo-{M}umford regularity of sheaves on surfaces}, American journal
  of mathematics,  (2006), pp.~373--417.

\bibitem{lant1}
\leavevmode\vrule height 2pt depth -1.6pt width 23pt, {\em Moduli spaces of
  principal bundles on singular varieties}, Kyoto Journal of Mathematics, 53
  (2013), pp.~3--23.


\bibitem{nagsesh1}
{\sc D.~S. Nagaraj and C.~S. Seshadri}, {\em Degenerations of the moduli spaces
  of vector bundles on curves {I}}, in Proceedings of the Indian Academy of
  Sciences-Mathematical Sciences, vol.~107, Springer, 1997, pp.~101--137.

\bibitem{ogmod}
{\sc K.~G. O'Grady}, {\em Moduli of vector bundles on projective surfaces: some
  basic results}, Inventiones mathematicae, 123 (1996), pp.~141--207.

\bibitem{piesch}
{\sc R.~Piene and M.~Schlessinger}, {\em On the {H}ilbert scheme
  compactification of the space of twisted cubics}, American Journal of
  Mathematics, 107 (1985), pp.~761--774.

\bibitem{sct4}
{\sc A.~Schmitt}, {\em On the modular interpretation of the
  {N}agaraj--{S}eshadri locus}, Journal f{\"u}r die reine und angewandte
  Mathematik (Crelles Journal), 2012 (2012), pp.~145--172.

\bibitem{S1}
{\sc E.~Sernesi}, {\em Deformaions of Algebraic Schemes}, Grundlehren der
  Mathematischen Wissenschaften-334, Springer-Verlag, 2006.

\bibitem{simp1}
{\sc C.~Simpson}, {\em Moduli of representations of the fundamental group of a
  smooth projective variety {I}}, Publications Math{\'e}matiques de l'IH{\'E}S,
  79 (1994), pp.~47--129.

\bibitem{simp2}
\leavevmode\vrule height 2pt depth -1.6pt width 23pt, {\em Moduli of
  representations of the fundamental group of a smooth projective variety
  {II}}, Publications Math{\'e}matiques de l'IH{\'E}S, 80 (1994), pp.~5--79.

\bibitem{bxia}
{\sc B.~Xia}, {\em Hilbert scheme of twisted cubics as a simple wall-crossing},
  Transactions of the American Mathematical Society, 370 (2018),
  pp.~5535--5559.

\end{thebibliography}
\end{document}